\documentclass[12pt]{article}
\usepackage{amsmath,amsfonts,amssymb,amsthm}
\usepackage[left=20mm,top=20mm,right=20mm,bottom=20mm,a4paper]{geometry}

\newcommand{\norm}[1]{\lVert#1\rVert}
\newcommand{\nnorm}[1]{\Bigl\|#1\Bigr\|}
\newcommand{\abs}[1]{\lvert#1\rvert}%

\def\rr{\mathbf{R}}
\def\qq{\mathbf{Q}}
\def\zz{\mathbf{Z}}

\def\F{\mathcal{F}}
\def\E{\mathrm{E}}
\def\angles#1{\langle{#1}\rangle}
\def\ffi{\varphi}
\DeclareMathOperator{\sgn}{sgn}

\theoremstyle{plain}  % italic
\newtheorem{theorem}{\hspace{\parindent}Theorem}%[section]
\newtheorem{lemma}[theorem]{\hspace{\parindent}Lemma}
\newtheorem{corollary}[theorem]{\hspace{\parindent}Corollary}

\theoremstyle{definition} % roman

\newtheorem{note}[theorem]{\hspace{\parindent}Note}

\renewenvironment{proof}[1][Proof]
  {\noindent\hspace{\parindent}\mbox{\textbf{#1.~}}}
  {\hfill$\triangleleft$} % ч.т.д. символ

\begin{document}
\hyphenation{}%

\title{
 Characterizing Hilbert spaces using Fourier transform over the field of
 p-adic numbers
 \footnote{
	   numbers.
 }
}
\author{Yauhen Radyna, Yakov Radyno, Anna Sidorik}
\date{}
\maketitle

\begin{abstract}
We characterize Hilbert spaces in the class of all Banach spaces
using Fourier transform of vector-valued functions over the field
$\qq_p$ of $p$-adic numbers. Precisely, Banach space $X$ is
isomorphic to a Hilbert one if and only if Fourier transform $F:
L_2(\qq_p,X)\to L_2(\qq_p,X)$ in space of functions, which are
square-integrable in Bochner sense and take value in $X$, is a
bounded operator.
\end{abstract}

In his paper~\cite{Peetre} J.~Peetre considered generalization of
Hausdorff--Young theorem  describing an image of  space $L_p(\rr)$
under Fourier transform. He considered  vector-valued functions
$x\in L_p(\rr,X)$, $1\le p\le 2$ over the real axis taking value
in Banach space $(X,\norm{\cdot})$ and integrable in Bochner
sense~\cite{Bochner-Integration}, i.e. weakly measurable and with
finite norm
\begin{equation}
 ||x||_{L_p(\rr,X)}=\Bigl(\int_{\rr}{||x(t)||^p d{t}}\Bigr)^{1/p}.
\end{equation}
For the case  $p=2$ J.~Peetre made the following observation. In
all known  cases Fourier transform
\begin{equation}
 \F:L_2(\rr,X)\to L_2(\rr,X), \quad
 (\F x)(s)= \int_{\rr}{x(t) e^{-2\pi i s t} d{t}}.
\end{equation}
was a bounded operator if and only if Banach space $X$ was
linearly isomorphic a Hilbert one. In~paper~\cite{Kwapien} Polish
mathematician S.~Kwapie\'n proofed that observation. Precisely,
one has

\begin{theorem}
The following statements are equivalent:

1) Banach space $X$ is isomorphic to a Hilbert one.

2) There exists $C>0$ such that for any positive integer  $n$ and
$x_0,x_1,x_{-1}, \ldots,$ $x_{n},x_{-n} \in X$
\begin{equation}
\int_0^1{\nnorm{\sum_{k=-n}^{n} e^{2\pi i k t}\cdot x_k}^2
d{t}}\le C\sum_{k = -n}^{n}\norm{x_k}^2.
\end{equation}

3) There exists $C>0$ such that for any positive integer $n$ and
$x_0,x_1,x_{-1}, \ldots,$ $x_{n},x_{-n} \in X$
\begin{equation}
\int_0^1{\nnorm{\sum_{k=-n}^{n} e^{2\pi i k t}\cdot x_k}^2
d{t}}\ge C^{-1}\sum_{k = -n}^{n}\norm{x_k}^2.
\end{equation}

4) Fourier transform $\F$, which is defined on a dense subset
$D_{\F}\subset L_2(\rr,X)$  consisting of elementary functions
\begin{equation}
 D_{\F}= \Bigl\{ x(t)= \sum_{k=1}^{n}{I_{A_k}(t)\cdot x_k}
 \Bigr\},
\end{equation}
where  $A_k$~are non-intersecting measurable subsets with finite
measure in $\rr$, $I_{A_k}$ are their indicators and $x_k\in X$,
is a bounded operator.
\end{theorem}

Now it's natural to look at Fourier transform on space of
vector-valued functions, which have their arguments not in $\rr$
but in a locally compact group $G$, and which are square
integrable in Bochner sense:
\begin{equation}
 \F= \F_G: L_2(G,X) \to L_2(\widehat{G},X),\quad
 (\F x)(s)=\int_{G}{x(t)\angles{s,t} d{\mu_G(t)}}.
\end{equation}
Here $X$ denotes Banach space, $\widehat{G}$ denotes Pontryagin
dual to $G$, $\angles{s,t}$ denotes canonical pairing between
$\widehat{G}$ and $G$, $\mu_G$ denotes Haar measure on $G$. It's
usefull to take the following dense subset
\begin{equation}
 D_\F= L_2(G)\otimes X \subset L_2(G,X),
\end{equation}
as a domain for $\F$ where its  action is described by formula
\begin{equation}
 \F\Bigl(\sum_{k=1}^{n}{\ffi_k(t)\cdot x_k} \Bigr)=
 \sum_{k=1}^{n}{((\F\ffi_k)(s)\cdot x_k)}.
\end{equation}

For a finite group $G$ a space $L_2(G,X)$ is isomorphic to a
finite direct product $X^G$ and Fourier transform is a bounded
operator for any Banach space $X$.

In this paper we consider the case of $G$ being a group $\qq_p$ of
$p$-adic numbers. Pontryagin dual group for $\qq_p$ can be
identified with $\qq_p$ itself using pairing
$\angles{s,t}=\chi_p(s t)= e^{2\pi i \{s t\}_p}$ where
$\{\cdot\}_p$~is a $p$-adic fractional part. Haar measure
$d{t}=d\mu(t)$ normed by $\mu(\zz_p)=1$ is
self-dual~\cite{VVZ,RadynaPadics}.

\begin{theorem}[Main]\label{th:main}
The following statements are equivalent:

1) Banach space  $X$ is linearly isomorphic to a Hilbert one.

2) There exists $C>0$ such that for any positive integer  $N$ and
$x_0,x_1, \ldots,$ $x_{p^{2 N}-1} \in X$
\begin{equation}
 \int_{\zz_p}
 {\nnorm{\sum_{k=0}^{p^{2 N}-1}\chi_p\Bigl(\frac{k t}{p^{2 N}}\Bigr)\cdot x_k}^2 d{t}}
 \le C\sum_{k = 0}^{p^{2 N}-1}\norm{x_k}^2.
\end{equation}

3) There exists $C>0$ such that for any positive integer  $N$ and
$x_0,x_1, \ldots,$ $x_{p^{2 N}-1} \in X$
\begin{equation}
 C^{-1}\sum_{k=0}^{p^{2 N}-1}\norm{x_k}^2 \le
 \int_{\zz_p}{\nnorm{\sum_{k=0}^{p^{2 N}-1}\chi_p\Bigl(\frac{k t}{p^{2 N}}\Bigr)\cdot x_k}^2 d{t}}.
\end{equation}

4) Fourier transform $\F: L_2(\qq_p,X)\to L_2(\qq_p,X)$ is a
bounded operator.
\end{theorem}

To proof Theorem~\ref{th:main} we use a criterion
from~\cite{Kwapien} formulated in probability terms.

\begin{lemma}\label{lm:Khinchin-type-ineq}

The following statements are equivalent:

1) Banach space  $X$ is linearly isomorphic to a Hilbert one.

2) There exists $C>0$ such that for any finite set of vectors
$x_1, x_2, ..., x_n\in X$ there holds a two-sided Khinchin-type
inequality
\begin{equation}
 C^{-1}\sum_{i=1}^{n}{\|x_i\|^2}\le \E\Bigl\|\sum_{i=1}^{n}{\delta_i
 x_i}\Bigr\|^2 \le C\sum_{i=1}^{n}{\|x_i\|^2}
\end{equation}
where  $\delta_i$ are independent random variables taking value
$\{+1,-1\}$ with probability~$\frac{1}{2}$, and symbol $\E$ stands
for expectation.
\end{lemma}

In order to use this criterion we construct a system of functions
on a set $\zz_p$ of $p$-adic integers which is analogous to
Rademacher system on  $[0,1]$
\begin{equation}
 r_i^\infty(z)= \sgn{\sin{2^i z}},\quad i=1,2,\ldots\;.
\end{equation}
For any $t\in\zz_p$ there exists canonical expansion  $t=t_0+
t_1\cdot p+ t_2\cdot p^2+\ldots$ where $t_k\in\{0,1,\ldots,p-1\}$.
Let's consider map
\begin{equation}
 \tau:\zz_p\to [0,1]: \sum_{k=0}^{+\infty}{t_k\cdot p^k}\mapsto
 \frac{1}{p}\sum_{k=0}^{+\infty}{t_k\cdot p^{-k}}.
\end{equation}
Restriction $\tilde\tau$ of this map onto  set
$\tilde\zz_p=\zz_p\setminus\{-1,-2,\ldots\}$ of full measure is a
bijection onto interval $[0,1)$. In addition $\tilde\tau$ is an
isomorphism of spaces with measure. Indeed, closed ball in
$\tilde\zz_p$ of measure $p^{-n}$ is nothing else but a set of
numbers with a fixed digits in canonical expansion starting from
$0$-th place to $(n-1)$-th one. The map $\tilde\tau$ sends this
set into a subset of reals with a fixed digits at places starting
from $0$-th to $(n-1)$-th one in base $p$ notation. The latter
subset also has measure $p^{-n}$. Balls in $\tilde\zz_p$ form a
base of topology,  thus $\tilde\tau$ preserves measure of all
Borel and Lebesgue-measurable sets, and $\tilde\tau$ is indeed an
isomorphism. So is $\tau$.

System of functions $r_i(t)=r_i^\infty(\tau(t))$ gives us a
realization of  $\delta_i$ of Lemma~\ref{lm:Khinchin-type-ineq}.
Thus
\begin{equation}
 \E\Bigl\|\sum_{i=1}^{n}{\delta_i\cdot x_i}\Bigr\|^2= \
 \int_{\zz_p}{\nnorm{\sum_{i=1}^{n}{r_i(t)\cdot x_i}}^2 d{t}}.
\end{equation}

Following~\cite{Kwapien} we proof a lemma for orthonormal systems
on $\zz_p$.

\begin{lemma}
Let $X$ be Banach space and let $(f_i)$ be a complete orthonormal
system in  $L_2(\zz_p)$. If for any $C>0$ and any vectors
$x_1,x_2, \ldots,x_n \in X$ there holds a two-sided inequality
\begin{equation}\label{eq:Rademacher-ortho-1}
C^{-1}\sum_{i=1}^{n}\norm{x_i}^2\le
\int_{\zz_p}{\nnorm{\sum_{i=1}^{n}f_i(t)\cdot x_i}^2 d{t}}\le
C\sum_{i=1}^{n}\norm{x_i}^2,
\end{equation}
then for the same constant $C>0$ and any vectors $x_1,x_2,
\ldots,x_n \in X$ we also have
\begin{equation}\label{eq:Rademacher-ortho-2}
C^{-1}\sum_{i=1}^{n}\norm{x_i}^2\le \int_{\zz_p}{\nnorm{\sum_{i =
1}^{n}r_i(t)\cdot x_i}^2 d{t}}\le C\sum_{i = 1}^{n}\norm{x_i}^2.
\end{equation}
Moreover, left inequality in~\eqref{eq:Rademacher-ortho-1} implies
left inequality in~\eqref{eq:Rademacher-ortho-2}, and right
inequality implies right one, independently of each other.
\end{lemma}

\begin{proof}
Orthonormality of $(r_i)$ and completeness of  $(f_k)$ imply that
for any  $\varepsilon>0$ we can find (using a standard gliding
hump procedure) a growing sequence of indexes $(k_j)$, $(m_j)$ and
orthogonal sequence $(h_j)$ such that
\begin{equation}
h_j=\sum_{k=k_j}^{k_{j+1}-1}(h_j,f_k)\cdot f_k,\
\int_{\zz_p}\abs{h_j(t)-r_{m_j}(t)}^2dt<\frac{\varepsilon}{2^{j}}.
\end{equation}
For any fixed  $x_1,x_2, \ldots,x_n \in X$ we have
\begin{equation}
\int_{\zz_p}{\nnorm{\sum_{i=1}^{n}r_i(t)\cdot x_i}^2 d{t}} =
\int_{\zz_p}{\nnorm{\sum_{i=1}^{n}r_{m_j}(t)\cdot x_i}^2 d{t}}.
\end{equation}
Triangle inequality yields
\begin{multline}
 \Bigl(\int_{\zz_p}\nnorm{\sum_{j=1}^{n}r_{m_j}(t)\cdot x_j}^2 d{t}\Bigr)^{1/2}\le
 \Bigl(\int_{\zz_p}\nnorm{\sum_{j=1}^{n}(r_{m_j}(t)-h_j(t))\cdot x_j}^2
 d{t}\Bigr)^{1/2} +
\\
 \Bigl(\int_{\zz_p}\nnorm{\sum_{j=1}^{n}h_{j}(t)\cdot x_j}^2 d{t}\Bigr)^{1/2}
 \le \sqrt{\varepsilon}\Bigl(\sum_{j=1}^{n}\norm{x_j}^2\Bigl)^{1/2}
 +\Bigl(\int_{\zz_p}\nnorm{\sum_{j=1}^{n}h_{j}(t)\cdot x_j}^2 d{t}\Bigr)^{1/2}.
\end{multline}
Equality $\displaystyle
1=\norm{h_j}^2=\sum_{k=k_j}^{k_{j+1}-1}\abs{(h_j,f_k)}^2$ yields
\begin{multline}
\int_{\zz_p}\nnorm{\sum_{j = 1}^{n}h_{j}(t)\cdot
x_j}^2dt=\int_{\zz_p}\nnorm{\sum_{j = 1}^{n}\biggl(
\sum_{k=k_j}^{k_{j+1}-1}(h_j,f_k)\cdot f_k\biggr)x_j}^2 d{t} \le
\\
 C \sum_{j = 1}^{n}\sum_{k=k_j}^{k_{j+1}-1}\abs{(h_j,f_k)}^2\norm{x_j}^2=
 C \sum_{j = 1}^{n}\norm{x_j}^2.
\end{multline}
Thus,
\begin{equation}
\int_{\zz_p}{\nnorm{\sum_{i = 1}^{n}r_i(t)\cdot x_i}^2 d{t}}\le
(\sqrt{\varepsilon}+\sqrt{C})^2\sum_{i = 1}^{n}\norm{x_i}^2.
\end{equation}
Since $\varepsilon$ is arbitrary small, we obtain the required
inequality. Right inequality in~\eqref{eq:Rademacher-ortho-2} can
be deduced in the same way.
\end{proof}

\begin{corollary}\label{cor:khinchin-ortho}
Let  $(f_i)$ be a complete orthonormal system in  $L_2(\zz_p)$ and
let  $X$ be a Banach space. Then space $X$ is linearly isomorphic
to a Hilbert one if and only if there exists constant $C>0$ such
that for any vectors $x_1,x_2, \ldots,x_n \in X$ the following
two-sided inequality holds:
\begin{equation}\label{eq:khinchin-r}
C^{-1}\sum_{i = 1}^{n}\norm{x_i}^2\le\int_{\zz_p}{\nnorm{\sum_{i =
1}^{n} f_i(t)\cdot x_i}^2 d{t}}\le C\sum_{i = 1}^{n}\norm{x_i}^2.
\end{equation}
\end{corollary}

\begin{lemma}\label{lm:F-calc}
Let
\begin{equation}\label{eq:F-calc}
 h(t)=\sum_{k=0}^{p^{2 N}-1}
 {p^{N/2} I_{B\left[\frac{k}{p^{N}},\frac{1}{p^{N}}\right]}(t)\cdot x_k},
 \quad x_k\in X.
\end{equation}

Then
\begin{equation}
 \norm{h}^2 = \sum_{k=0}^{p^{2 N}-1}\norm{x_k}^2, \quad
 \norm{\F h}^2 =
 \int_{\zz_p}{\nnorm{\sum_{k=0}^{p^{2 N}-1}\chi_p\Bigl(\frac{k t}{p^{2 N}}\Bigr)\cdot x_k}^2 d{t}}.
\end{equation}
\end{lemma}

\begin{proof}
By direct calculations one has
\begin{equation}
 \norm{h}^2=
 \sum_{k=0}^{p^{2 N}-1}{p^N \mu\Bigl(B\Bigl[\frac{k}{p^N},\frac{1}{p^N}\Bigr]\Bigr)\cdot\norm{x_k}^2}=
 \sum_{k=0}^{p^{2 N}-1}\norm{x_k}^2,
\end{equation}
\begin{multline}
 \nnorm{\F h}^2=
 \int_{\qq_p}
 {\nnorm{\sum_{k=0}^{p^{2 N}-1}{p^{N/2} \Bigl(\F I_{B\left[\frac{k}{p^N},\frac{1}{p^N}\right]}\Bigr)(t)
 \cdot x_k}}^2 d{t}} =
 \\
 \int_{\qq_p}
 {\nnorm{\sum_{k=0}^{p^{2 N}-1}
 {p^{-N/2} I_{B\left[0,p^N\right]}(t)\cdot\chi_p\Bigl(\frac{k t}{p^N}\Bigr)\cdot x_k}}^2 d{t}}
 = \int_{\zz_p}{\nnorm{\sum_{k=0}^{p^{2 N}-1}\chi_p\Bigl(\frac{k t}{p^{2 N}}\Bigr)\cdot x_k}^2 d{t}}.
\end{multline}
\end{proof}

\begin{lemma}\label{lem:ineq-dual}
Let $X$ be Banach space, let $(f_i)$ be complete orthonormal
system in $L_2(\zz_p)$, and let $C>0$ be a constant such that for
any vectors $x_1,x_2, \ldots,x_n \in X$
\begin{equation}
 \int_{\zz_p}{\nnorm{\sum_{i=1}^{n}f_i(t)\cdot x_i}_X^2 d{t}}\ge
 C^{-1}\sum_{i=1}^{n}\norm{x_i}^2.
\end{equation}
Then for any vectors $x_1^*,x_2^*, \ldots,x_m^* \in X^*$ in dual
vector space we have reverse type inequality
\begin{equation}
 \int_{\zz_p}{\nnorm{\sum_{j=1}^{m}f_j(t)\cdot x_j^*}_{X^*}^2 d{t}}\le
 C\sum_{j=1}^{m}\norm{x_j^*}^2.
\end{equation}
\end{lemma}

\begin{proof}
Completeness of  $(f_i)$ implies that linear combinations of $f_i$
form a dense subset in $L_2(\zz_p)$. That's why set $E$ of all
vector-valued functions
\begin{equation} \varphi(t) =
\sum_{i=1}^{n}f_i(t)\cdot x_i,\quad x_i \in X,
\end{equation}
is dense in space $L_2(\zz_p,X)$ of Bochner square-integrable
functions over $\zz_p$. Definition of norm in the dual space
$L_2(\zz_p,X^*)=L_2(\zz_p,X)^*$ yields that for any fixed linear
combination $\displaystyle\varphi^*(t) =
\sum_{j=1}^{m}{f_j(t)\cdot x_j^*}$ and any $\varepsilon > 0$ there
exists such a $\displaystyle\varphi(t) =
\sum_{i=1}^{n}{f_i(t)\cdot x_i} \in E$ with unit norm
$\displaystyle ||\varphi||_{L_2(\zz_p,X)}=
\Bigl(\int_{\zz_p}{\norm{\varphi(t)}_X^2 d{t}}\Bigr)^{1/2}= 1$
that
\begin{multline}
 \Bigl(\int_{\zz_p}{\nnorm{\sum_{j=1}^{m}{f_j(t)\cdot x_j^*}}_{X^*}^2 d{t}}\Bigr)^{1/2}
 \le\int_{\zz_p}{\abs{\angles{\varphi^*(t),\varphi(t)}}d{t}}+\varepsilon
 =\sum_{i=1}^{\min(m,n)}{\abs{\angles{x_i^*,x_i}}}+\varepsilon \le
\\
 \le \Bigl(\sum_{j=1}^{m}{\norm{x_j^*}_{X^*}^2}\Bigr)^{1/2}
 \Bigl(\sum_{i=1}^{n}\norm{x_i}_X^2\Bigr)^{1/2}+\varepsilon
 \le{\sqrt{C}}\Bigl(\sum_{j=1}^{m}{\norm{x_j^*}_{X^*}^2}\Bigr)^{1/2}+\varepsilon.
\end{multline}
And we obtain the desired inequality.
\end{proof}

\bigskip
\begin{proof}[Proof~of~Theorem~\ref{th:main}]
Straightforward computation shows that if $X$ is linearly
isomorphic to a Hilbert space then for any orthonormal system of
functions  $f_1,f_2, \ldots f_n \in L_2(\zz_p)$ and vectors
$x_1,x_2,\ldots,x_n \in X$ there holds equality
\begin{equation}
\int_{\zz_p}{\nnorm{\sum_{i=1}^{n}{f_i(t)\cdot x_i}}^2 d{t}}
=\sum_{i=1}^{n}\norm{x_i}^2.
\end{equation}
System of characters $\chi_p(k t/p^{2 N})$, $k=0,1,\ldots,p^{2
N}-1$ is complete and orthonormal in $L_2(\zz_p)$, thus we have
implications $\mathit{1)}\Rightarrow \mathit{2)}$ и
$\mathit{1)}\Rightarrow \mathit{3)}$.

Corollary~\ref{cor:khinchin-ortho} yields an implication
$\mathit{2)}\;\&\;\mathit{3)} \Rightarrow 1)$.

Let us assume that condition $\mathit{2)}$ holds.
Lemma~\ref{lm:F-calc} shows that that for any function $h$ of
type~\eqref{eq:F-calc} there is an equality
\begin{equation}
 \norm{\F h}^2=
 \int_{\zz_p}{\nnorm{\sum_{k=0}^{p^{2 N}-1}{\chi_p\Bigl(\frac{k t}{p^{2 N}}\Bigr)\cdot x_k}}^2 d{t}}
 \le C \sum_{k=0}^{p^{2 N}-1}{\norm{x_k}^2}= C\norm{h}^2.
\end{equation}
Functions of type~\eqref{eq:F-calc} are elements of
$L_2(\qq_p)\otimes X$ and form a dense subset in $L_2(\qq_p,X)$.
Thus Fourier transform  $\F: L_2(\qq_p,X)\to L_2(\qq_p,X)$ is
bounded, and  $\mathit{2)} \Rightarrow \mathit{4)}$.

Let us assume that $\mathit{4)}$ holds, i.e. Fourier transform is
bounded. Than restriction of inverse Fourier transform $\F^{-1}=
\F^3$ onto a dense subset $L_2(\qq_p)\otimes X$ is also bounded,
and there exists a constant $C_1$ such that for any  $h\in
L_2(\qq_p)\otimes X$ there holds inequality $\norm{\F^{-1} h} \le
C_1 \norm{h}$. Let $h$ be of type~\eqref{eq:F-calc}. One has
\begin{equation}
 \sum_{k=0}^{p^{2 N}-1}{\norm{x_k}}^2=\norm{h}=\norm{\F^{-1}\F h}\le
 C_1 \norm{\F h}=
 C_1 \int_{\zz_p}{\nnorm{\sum_{k=0}^{p^{2 N}-1}{\chi_p\Bigl(\frac{k t}{p^{2 N}}\Bigr)\cdot x_k}}^2
 d{t}}.
\end{equation}
Thus $\mathit{4)}\Rightarrow \mathit{3)}$.

Let us assume that condition $\mathit{3)}$ holds.
Lemma~\ref{lem:ineq-dual} yields that the dual space $X^*$
satisfies a condition analogous to $\mathit{2)}$. Then there holds
a condition, which is analogous to $4)$, i.e. Fourier transform in
$L_2(\qq_p,X^*)$ is bounded. That is why $X^*$ satisfies a
condition analogous to~$\mathit{3)}$. Applying
Lemma~\ref{lem:ineq-dual} once more we see that $X^{**}$ satisfies
a condition analogous to~$\mathit{2)}$. So does $X\subset X^{**}$.
Thus  $\mathit{3)}\Rightarrow \mathit{2)}$.

Now we have enough implications to state that conditions
$\mathit{1)}$, $\mathit{2)}$, $\mathit{3)}$, $\mathit{4)}$ are all
equivalent.
\end{proof}

\begin{note}
In fact, condition $\mathit{2)}$ in Theorem~\ref{th:main} is
equivalent to boundedness of Fourier transform $\F_{\qq_p/\zz_p}:
L_2(\qq_p/\zz_p,X)\to L_2(\zz_p,X)$. Applying  Fourier transform
to a function with finite support in $L_2(\qq_p/\zz_p,X)$ we
directly obtain $\mathit{2)}$ from boundedness of $\F$. Since
finitely-supported functions form a dense subset in
$L_2(\qq_p/\zz_p,X)$, we have a reverse implication.

Inverse Fourier transform $\F_{\zz_p}^{-1}: L_2(\qq_p/\zz_p,X)\to
L_2(\zz_p,X)$ equals to composition $\F_{\zz_p}^{-1}=
\mathcal{I}_{\zz_p} \F_{\qq_p/\zz_p}$ where $\mathcal{I}_{\zz_p}$
is an isometry $x(t)\mapsto x(-t)$. Thus, it's also bounded.

Using locally-constant functions on $\zz_p$ we see that
$\mathit{3)}$ is equivalent to boundedness of Fourier transform
$\F_{\zz_p}: L_2(\zz_p,X)\to L_2(\qq_p/\zz_p,X)$ and inverse
Fourier transform $\F_{\qq_p/\zz_p}^{-1}: L_2(\zz_p,X)\to
L_2(\qq_p/\zz_p,X)$.
\end{note}

\def\refname{

%
%\begin{flushleft}\bigskip
%\textit{Mechanics and Mathemathics Faculty,\\
%Belarusian State University,\\
%F. Skaryna av., 4,\\
%BY-220030, Minsk, BELARUS.}
%\end{flushleft}
%
\end{document}